\documentclass[12pt]{article}
\begin{document}
\newcommand{\mm}{{\rm Mod}_S}
\title{Fifteen Problems about the Mapping Class Groups}
\date{}
\author{Nikolai V. Ivanov}

\renewcommand{\thefootnote}{\fnsymbol{footnote}}

\maketitle

\footnotetext{Supported in part by the NSF Grants DMS-9401284 and DMS-0406946.}

Let $S$ be a compact orientable surface, possibly with boundary.  We
denote by $\mm$ the {\em mapping class group} of $S$, i.e. the group
$\pi _0 ({\rm Diff} (S))$ of isotopy classes of diffeomorphisms of $S$.
This group is also known as the {\em Teichm\"uller modular group\/} of $S$,
whence the notation $\mm$. Note that we include the isotopy classes of
orientation-reversing diffeomorphisms in $\mm$, so our group $\mm$ is
what is sometimes called the {\em extended} mapping class group.

For any property of discrete groups one may ask if $\mm$ has this
property. Guidance is provided by the analogy, well-established by now,
between the mapping class groups and arithmetic groups. See,
for example, \cite{i-uspekhi} or \cite{i-handbook} for a discussion.  

After the 1993 Georgia Topology Conference, R. Kirby was preparing a new
version of his famous problem list in low-dimensional topology. In
response to his appeal, I prepared a list of ten problems about the
mapping class groups, which was informally circulated as a preprint \cite{i-problems}.
Some of these problems ended up in Kirby's list \cite{kirby} in a
somewhat modified form, some did not. In the present article I will indicate
the current status of these problems and also will present 
some additional problems. Only Problems 4,6, and 9 of the original list of ten problems 
are by now completely solved. (For the convenience of the readers familiar with \cite{i-problems} I preserved the numbering of these ten problems.)

I tried to single out some specific questions, leaving aside such well-known problems as
the existence of finite dimensional faithful linear representations
or the computation of the cohomology groups.

\paragraph{1. The Congruence Subgroups Problem.} This is Problem 2.10 in Kirby's
list.\\ 

Suppose that $S$ is closed. Recall that a subgroup $\Gamma$ of a group $G$ is called {\em
characteristic\/} if $\Gamma$ is invariant under all automorphisms of
$G$.  If $\Gamma$ is a characteristic subgroup of $\pi _1 (S)$, then
there is a natural homomorphism ${\rm Out}(\pi _1 (S))\rightarrow {\rm
Out} (\pi _1 (S) /\Gamma )$, where for a group $G$ we denote by ${\rm
Out} (G)$ the quotient of the group ${\rm Aut} (G)$ of all
automorphisms of $G$ by the (automatically normal) subgroup of all
inner automorphisms. Clearly, if $\Gamma$ is of finite index in $\pi _1
(S)$, then the kernel of this homomorphism is also of finite index.
Note that any subgroup of finite index in a finitely generated group
(in particular, in $\pi _1 (S)$) contains a characteristic subgroup of
finite index. Since by the Dehn--Nielsen theorem $\mm$ is canonically
isomorphic to ${\rm Out}(\pi _1 (S))$, our construction gives rise to a
family of subgroups of finite index in $\mm$.  By analogy with the
classical arithmetic groups we call them the {\em congruence
subgroups}. 

\paragraph{Conjecture.} {\em Every subgroup of finite index in $\mm$ contains a congruence subgroup.}\\

V. Voevodsky had indicated (in a personal communication) a beautiful
application of this conjecture. Namely, the conjecture implies that a
smooth algebraic curve over ${\bf Q}$ is determined up to
isomorphism by its algebraic fundamental group (which is isomorphic to
the profinite completion of $\pi _1 (S)$) considered together with the
natural action of the absolute Galois group ${\rm Gal} (\overline{\bf
Q} /{\bf Q})$ on it. This corollary was apparently first conjectured by A. Grothendieck. 
I am not aware of any publication where this conjecture of Grothendieck is deduced from the solution of the congruence subgroup problem for the mapping class groups.
The Grothendieck conjecture itself was proved by S. Mochizuki \cite{mo1} after the initial work by A. Tamagawa \cite{ta}. See also \cite{mo3} and the expository accounts by S. Mochizuki \cite{mo2} and G. Faltings \cite{faltings}. The proof of the Grothendieck conjecture lends some additional credibility (beyond the analogy with the classical congruence-subgroup problem) to the above 
congruence subgroup conjecture.

\paragraph{2. Normal Subgroups.} This is Problem 2.12 (B) in Kirby's list.\\

If a subgroup $\Gamma$ of $\pi _1 (S)$
is characteristic, then the kernel of the natural homomorphism ${\rm
Out}(\pi _1 (S))\rightarrow {\rm Out} (\pi _1 (S) /\Gamma )$ from the
Problem 1 is a normal subgroup of ${\rm Out}(\pi _1 (S))$. So, this
construction gives rise to a family of normal subgroups of $\mm$. In
general, these subgroups have infinite index. For example, the Torelli
subgroup is a subgroup of this type.\\ 

\paragraph{Question.} {\em Is it true that any normal subgroup
is commensurable with such a subgroup?}\\

Recall that two subgroups
$\Gamma _1$ and $\Gamma _2 $ of a group $G$ are {\em commensurable} if the
intersection $\Gamma _1 \cap \Gamma _2$ has finite index in both
$\Gamma_1$ and $\Gamma_2$.

This problem was suggested by a discussion with H. Bass.

\paragraph{3. Normal Subgroups and pseudo-Anosov Elements.} This is Problem
2.12 (A) in Kirby's list.\\  

\paragraph{Question.} {\em Is it possible that all nontrivial elements of a {\em normal\/} subgroup of $\mm$ are
pseudo-Anosov?}\\

To the best of my knowledge, this question was posed independently by D. D. Long, J. D. McCarthy, and R. C. Penner in the early eighties. I learned it from R. C. Penner in 1984.

Moderate progress in the direction of the solution of this problem is due to K. Whittlesey \cite{wh}, who constructed examples of such subgroups in the case of spheres with 5 or more punctures and for closed surfaces of genus 2 (in the latter case the mapping class group
is well-known to be intimately connected with the mapping class group of
the sphere with 6 punctures). Unfortunately, there is, apparently, no hope of extending the method of proof to the other cases, especially to the closed surfaces of higher genus (which were the
focus of all previous attempts at this problem, I suspect).

\paragraph{4. Conjecture (Mostow-Margulis superrigidity).} This is Problem 2.15
in Kirby's list.

\paragraph{Conjecture.}{\em If $\Gamma$ is an irreducible arithmetic group of rank $\geq 2$, then every
homomorphism $\Gamma \rightarrow \mm$ has finite image.}\\

For many arithmetic groups $\Gamma$ the conjecture can be proved by
combining some well-known information about arithmetic groups with
equally well-known properties of $\mm$. But, for example, for cocompact
lattices in $ SU(p,q)$ this straightforward approach seems to fail. (I
owe this specific example to G. Prasad.) In any case, by now this conjecture
is completely proved.

V. A. Kaimanovich and H. Masur \cite {kam}  proved that so-called non-elementary subgroups
(subgroups containing a pair of pseudo-Anosov elements with disjoint sets of fixed points
in the Thurston boundary) of $\mm$ are
not isomorphic to irreducible arithmetic groups of rank $\geq 2$. The
proof is based on the theory of random walks on $\mm$ developed by Masur
and Kaimanovich--Masur and on the results of H. Furstenberg about
random walks on arithmetic groups (see the references in \cite{kam}). If combined with the Margulis finiteness theorem and the techniques of \cite{i-book}, this result can be used to prove the conjecture.

A proof close in the spirit to the one alluded to in the previous paragraph is due to B. Farb and H. Masur \cite{fam}.

In addition, S-K. Yeung \cite{yeung} completed the picture by proving that any homomorphism from a lattice $\Gamma$ in either ${{\rm Sp}}(m,1)$ or $F_4^{-20}$ (the isometry group of the Cayley plane) into a mapping class group has finite image. The lattices in the remaining simple Lie groups of rank 1
often admit non-trivial homomorphisms to ${\bf Z}$, and, therefore, homomorphisms to the mapping class groups with infinite image.

\paragraph{5. Dehn multitwists.} For a nontrivial circle $\alpha$ on
$S$ (i.e. a submanifold $\alpha$ of $S$ diffeomorphic to the standard circle $S^1$) let us denote by $t_ \alpha$ the (left) Dehn twist along $\alpha$;
$t_ \alpha \in \mm$.  A {\em Dehn multi-twist} is any composition of
the form $t_{\alpha _1}^{\pm 1} \circ t_{\alpha _2}^{\pm 1} \circ
\cdots \circ t_{\alpha _n}^{\pm 1}$ for disjoint circles  $\alpha _1
,\alpha _2, \ldots ,\alpha _n$.

\paragraph{Question.} {\em Is there a constant $N_S$, depending only on
$S$, such that the following holds? Let $f \in \mm$ and let $t= t_{\alpha
_1}^{\pm 1} \circ t_{\alpha _2}^{\pm 1} \circ \cdots \circ t_{\alpha
_n}^{\pm 1}$ be a Dehn multi-twist. If $A=\{ f^m (\alpha _i) : 1 \leq i
\leq n , m\in {\bf Z} \}$ fills $S$ (i.e. for any nontrivial circle
$\gamma$ there exists an $\alpha \in A$ such that the geometric intersection number ${\rm i} (\gamma
,\alpha ) \neq 0 $), then only a finite number of elements of the form $t^j\circ
f$, $j\in {\bf Z}$, are not pseudo-Anosov, and the corresponding $j$ are among some $N_S$
consecutive integers.}\\

By a theorem of A. Fathi \cite{fathi}, if $t$ is a Dehn {\em twist}, then this is true for $N_S =7$. This theorem of A. Fathi was preceded by a theorem of D. Long and H. Morton \cite{longm} to the effect that for a Dehn twist $t$, only a finite number of element $t^j\circ f$ are not pseudo-Anosov under the above assumptions, without a uniform bound on the number of exceptions.\\

A weaker version of this question is still interesting: under the same
conditions, is it true that no more than $N_S$ of $t^j
\circ f$ are not pseudo-Anosov? 

The initial motivation for this question was that a positive answer
would allow to prove the Conjecture 4 in some nontrivial cases. This motivation
is now (completely?) obsolete, but I still consider this question to be interesting. I believe that
Fathi's paper \cite{fathi} is one of the deepest and the most underappreciated works
in the theory of the mapping class groups. In fact, I am aware of only one application 
of his results: the author's theorem \cite{i-rank} to the effect that the mapping class groups have rank 1 in an appropriate sense (note that the Long--Morton theorem \cite{long} is not
sufficient here); see \cite{i-handbook}, Section 9.4 for a discussion. The above question may be considered as a test of our understanding of Fathi's ideas.

\paragraph{6. Automorphisms of Complexes of Curves.} The {\em complex
of curves} $C(S)$ of a surface $S$ is a simplicial complex defined as
follows. The vertices of $C(S)$ are the isotopy classes of nontrivial
(i.e., not bounding a disk and non deformable into the boundary)
circles on $S$. A set of vertices forms a simplex if and only if they
can be represented by disjoint circles. This notion was introduced by
W. J. Harvey \cite{harv}. Clearly, $\mm$ acts on $C(S)$, and
this action is almost always effective (the main exception is the case
of a closed surface of genus 2, where the hyperelliptic involution acts on $C(S)$
trivially). If the genus of $S$ is at least 2, then all automorphisms of $C(S)$ come
from $\mm$, according to a well-known theorem of the author \cite{i-ihes},
\cite{i-imrn}. The problem was to prove the same for surfaces of genus 1 and 0.

In his 1996 MSU Ph.D. thesis M. Korkmaz proved that all
automorphisms of $C(S)$ come from $\mm$ for all surfaces of genus $0$
and $1$ with the exception of spheres with $\leq 4$ holes and tori with
$\leq 2$ holes. See his paper \cite{kork}.

Later on, F. Luo \cite{luo} suggested another
approach to the above results about automorphisms of $C(S)$, still based on the ideas of
\cite{i-ihes}, and also on a multiplicative structure
on the set of vertices of $C(S)$ introduced in \cite{luo1}.
He also observed 
that ${\rm Aut ( C(S) )}$ is not equal to $\mm$ if $S$ is a torus with $2$
holes. The reason is very simple: if $S_{1,2}$ is a torus with $2$
holes, and $S_{0,5}$ is a sphere with $5$ holes, then $C(S_{1,2})$ is
isomorphic to $C(S_{0,5})$, but ${\rm Mod} _{S_{1,2}}$ is not isomorphic to
${\rm Mod} _{S_{0,5}}$. 

So, the problem is solved completely. The theorem about the automorphisms of $C(S)$ had stimulated a lot of further results about the automorphisms of various objects 
related to surfaces; the proofs are usually based on a reduction of the problem in question to the theorem about the automorphisms of $C(S)$. I will not attempt to survey these results. Instead of this, I will state a metaconjecture.\\

\paragraph{Metaconjecture.} {\em Every object naturally associated to a surface $S$ and having a sufficiently rich structure has $\mm$ as its groups of automorphisms. Moreover, this can be
proved by a reduction to the theorem about the automorphisms of $C(S)$.}

\paragraph{7. The First Cohomology Group of Subgroups of Finite Index.} This is Problem 2.11 (A) in Kirby's list.\\

\paragraph{Question.}{\em Is it true that $H^1 (\Gamma)=0$ for any subgroup $\Gamma$ of finite
index in $\mm$?}\\

It is well known that $H^1 (\mm )=0$. In his 1999 MSU Ph.D. thesis 
F. Taherkhani carried out some extensive computer 
calculation aimed at finding a subgroup of finite index with non-zero first cohomology group.
For genus 2 he found several subgroups $\Gamma$ with $H^1 (\Gamma) \neq 0$, but in genus 3 all examined subgroups $\Gamma$ turned out to have $H^1 (\Gamma)=0$. The higher genus cases
apparently were well beyond the computer resources available at the time. See
\cite{tah}.

J. D. McCarthy \cite{mc} proved that if $S$ is a closed surface of genus $\geq 3$
and $\Gamma$ is a subgroup of finite index in $\mm$ containing the
Torelli subgroup, then the first cohomology group $H^1 (\Gamma )$ is
trivial. His methods are based on D. Johnson results about the Torelli
subgroup, the solution of the congruence subgroups problem for $Sp_{2g}
({\bf Z})$, $g\geq 3$, and the Kazhdan property (T) of $Sp_{2g} ({\bf
Z})$, $g\geq 3$.

\paragraph{8. Kazhdan Property (T).} This is Problem 2.11 (B) in Kirby's list.

\paragraph{Question.}{\em Does $\mm$ have the Kazhdan property (T)?}\\ 

A positive answer would imply the positive answer to the previous question,
but this problem seems to be much more difficult. Most of the known
proofs of the property (T) for discrete groups are eventually based on
the relations of these discrete groups with Lie groups and on the representation theory of Lie
groups. Such an approach is not available for $\mm$. New approaches to the
Kazhdan Property (T) (see, for example Y. Shalom \cite{sha}, A. \.Zuk \cite{zuk}
and the N. Bourbaki Seminar report of A. Valette \cite{val}) hold a better
promise for our problem, but, to the best of my knowledge, no serious work
in this direction has been done.

The problem remains completely open.

\paragraph{9. Unipotent Elements.} This is Problem 2.16 in Kirby's list. 

Let $d_W (\cdot ,\cdot )$ be the
word metric on $\mm$ with respect to some finite set of generators. Let
$t\in \mm $ be a Dehn twist. 

\paragraph{Question.}{\em What is the growth rate of $d_W (t^n , 1)
$?}\\

One would expect that either the growth is linear, or $d_W (t^n ,
1)=O(\log n) $. In the arithmetic groups case, the logarithmic growth
corresponds to virtually unipotent elements of arithmetic groups of
rank $\geq 2$, according to a theorem of A. Lubotzky, S. Moses
and M.S. Raghunathan \cite{lmr}.

B. Farb, A. Lubotzky and Y. Minsky \cite{flm} proved that the growth
is linear, so the problem is solved completely. 

\paragraph{10. Nonorientable Surfaces.} Most of the work on the mapping
class groups is done for orientable surfaces only. Some exceptions are
the paper of M. Scharlemann \cite{shar}, and the computation of the virtual cohomological dimension for the mapping class groups of nonorientable surfaces in the author's work \cite{i-uspekhi}. One may look for analogues of other theorems about $\mm$ for nonorientable surfaces. For example,
what are the automorphisms of mapping class groups of nonorientable
surfaces?

Of course, one may expect that some results will extend to the nonorientable case more or
less automatically. This applies, for example, to the author's results about the subgroups of
the mapping class groups reported in \cite{i-book}. The only reason why these results are stated only in the orientable case is the fact that expositions of Thurston's theory of surfaces exist only in the orientable case. But in other problems some new phenomena appear. See, for example, the result of M. Korkmaz \cite{kork2}, \cite{kork3}. N. Wahl \cite{whal} recently reported some work in progress aimed at proving a homology stability theorem for the mapping class groups of nonorientable surfaces, where some new difficulties appear in comparison with the 
orientable case. 

\paragraph{11. Free subgroups generated by pseudo-Anosov elements.} One of the first results
in the modern theory of the mapping class groups was the theorem to the effect that for two independent pseudo-Anosov elements (where {\em independence} means that their sets of fixed points in the Thurston boundary are disjoint) their sufficiently high powers generate a free group (with two generators). This observation was done independently by J.D. McCarthy and the author and for the author it was the starting point of the theory presented in \cite{i-book}. This result clearly extends to any finite number of independent elements and even to an infinite collection of pairwise independent elements, {\em provided that we do not impose a bound} on the powers to which our elements are raised. Now, what if we do impose a uniform bound?\\

\paragraph{Conjecture.}{\em If $f$ is pseudo-Anosov element of a mapping class group $\mm$
with sufficiently large dilatation coefficient, then the subgroup of $\mm$ normally generated by $f$ is a free group having as generators the
conjugates of $f$. More cautiously, one may conjecture that the above holds for a sufficiently high power $g=f^N$ of a given pseudo-Anosov element $f$.}\\

This conjecture is motivated by a theorem of M. Gromov (see \cite{g}, Theorem 5.3.E).
According to this theorem, the subgroup of a hyperbolic group normally generated by a 
hyperbolic element (with sufficiently big translation length) is a free group having as generators the conjugates of this element. Of course, it is well-known that mapping class groups are not hyperbolic. But, as M. Gromov noticed, what is essential for his proof is not the
global negative curvature (hyperbolicity), but the negative curvature
around the loop representing the considered element in an Eilenberg-MacLane space of the group in question. So, the hope is that the Teichm\"uller spaces have enough negative curvature
along the axes of pseudo-Anosov elements for a similar conclusion to hold. For a more direct
and elementary approach to Gromov's result see the work of Th. Delzant \cite{d}.

This conjecture may be relevant to Problem 3 above, since one may expect that all 
nontrivial elements of such subgroups are pseudo-Anosov.

\paragraph{12. Deep relations.} It is well-known that mapping class groups of closed surfaces are generated
by the Dehn twists (this result is due to Dehn himself, see \cite{i-handbook} for a detailed
discussion), and that some specific relations among Dehn twists, such as the Artin (braid) relations and the lantern relations play a crucial role 
in the proofs of many results about the mapping
class groups. In fact, according to a theorem of S. Gervais \cite{ge} (a hint to such a theorem is contained already in J. Harer's paper \cite{ha}) $\mm$ admits a presentation having all Dehn twists as generators and only some standard relations between them as relations. There is also a similar presentation having only Dehn twists about nonseparating curves as generators. See \cite{ge} for exact statements. (For surfaces with boundary simple additional generators are needed, but the following discussion of the relations between Dehn twists equally makes sense for surfaces with boundary.)

Almost all of these relations disappear if we replace Dehn twists by powers of Dehn twists (for example,
if we are forced to do so by trying to prove something about subgroups of finite index). Apparently, only the relations 
\[T_{\alpha} T_{\beta} T_{\alpha}^{-1}=T_{T_{\alpha}({\beta})},\]
survive, where $T_{\gamma}$ for a circle $\gamma$ denotes some power $t_{\gamma}^L$ of a Dehn twist $t_{\gamma}$ about ${\gamma}$. We are primarily interested in the case of powers independent of circles, since there are no other natural way to assign powers to circles. Notice that these relations include the commutation relations $T_{\alpha} T_{\beta}=T_{\beta} T_{\alpha}$ for disjoint ${\alpha},{\beta}$, and for the circles ${\alpha},{\beta}$ intersecting once transversely they imply the Artin relation $t_{\alpha} t_{\beta} t_{\alpha} = t_{\beta} t_{\alpha} t_{\beta}$.

\paragraph{Question.} {\em Are there any other relations between $N$-th powers $T_{\gamma}=t_{\gamma}^N$  of Dehn
twists for sufficiently large $N$? In other words, do the above relations provide a presentation
of the group generated by the $N$-th powers of Dehn twists? If not, what are the additional relations?}\\

\paragraph{13. Burnside groups.} Note that the subgroup generated by the $N$-th powers of Dehn twists is often of infinite index in $\mm$, as it follows from the results of L. Funar \cite{fun}. (See also G. Masbaum \cite{mas} for a simple proof of more precise results.)\\

\paragraph{Question.}{\em Is the subgroup of $\mm$ generated by the $N$-th powers of {\em all} elements of $\mm$ of infinite index in $\mm$ for sufficiently large $N$?}\\

Notice that such a subgroup is obviously normal and if it is of infinite index, then the quotient group is an infinite Burnside group.

\paragraph{14. Homomorphisms between subgroups of finite index.} In general, one would like to understand homomorphisms between various mapping class groups and their finite index subgroups. Here is one specific conjecture.

\paragraph{Conjecture.} {\em Let $S$ and $R$ be {\em closed} surfaces. Let $\Gamma$ be a subgroup of finite index in $\mm$. If the genus of $R$ is less than the genus of $S$, then there is no homomorphism $\Gamma \rightarrow {{\rm Mod}_R}$ having as an image a subgroup of finite index in
${{\rm Mod}_R}$.}\\

In fact, one may hope that any such homomorphism $\Gamma \rightarrow {{\rm Mod}_R}$ has a finite image. But while Problem 7 about $H^1(\Gamma)$ is unresolved, a more cautious conjecture seems to be more appropriate and more accessible, because if $H^1(\Gamma)$ is infinite, $\Gamma$ admits a homomorphism onto $\bf Z$, and therefore, a lot of homomorphisms into ${{\rm Mod}_R}$, in particular, with infinite image.

In the special case when $\Gamma=\mm$, the conjecture was recently proved by W. Harvey and M. Korkmaz \cite{harvk}. Their methods rely on the use of elements of finite order in $\Gamma=\mm$, and therefore cannot be extended to general subgroups $\Gamma$ of finite index.

\paragraph{15. Frattini subgroups.} Recall that the Frattini subgroup $\Phi (G)$ of a group $G$ is the intersection of all proper (i.e. different from $G$) maximal subgroups of $G$. According to a theorem of Platonov \cite{platonov} the Frattini subgroup $\Phi (G)$ of a finitely generated linear group $G$ is nilpotent. An analogue of this result for subgroups of mapping class groups was proved by the author \cite{i-book} after some initial results of D. D. Long \cite{long}. Namely, for any subgroup $G$ of $\mm$, its Frattini subgroup $\Phi (G)$ is nilpotent; see \cite{i-book}, Chapter 10.

Despite the similarly sounding statements, the theorems of Platonov and of the author have a completely different nature. In order to explain this, let us denote by $\Phi _{f} (G)$ the intersection of all maximal subgroups of {\em finite index} in $G$, and by $\Phi _{i} (G)$ the intersection of all maximal subgroups of {\em infinite index} in $G$; clearly, $\Phi (G)= \Phi _{f} (G)\cap \Phi _{i} (G)$. It turns out that the main part of the argument in \cite{i-book} proves that $\Phi _{i} (G)$ is virtually abelian for subgroups of $\mm$ (and so $\Phi (G)$ is both virtually abelian and nilpotent). At the same time, Platonov actually proves that $\Phi _{f} (G)$ is nilpotent for finitely generated linear groups $G$; subgroups of infinite index play no role in his arguments. (Actually, the very existence of maximal subgroups of infinite index in linear groups is a highly nontrivial result due to G. A. Margulis and G. A. Soifer \cite{mar-soifer}, and it was proved much later than Platonov's theorem.) 

\paragraph{Conjecture.} {\em For every finitely generated subgroups $G$ of $\mm$, the group 
$\Phi _{f} (G)$ is nilpotent.}\\

For a more detailed discussion, see \cite{i-book}, Section 10.10. I repeated here
this old problem in order to stress that our understanding of the subgroups of finite index in
$\mm$ is rather limited.

\noindent
{\sc Michigan State University\\
Department of Mathematics\\
Wells Hall\\
East Lansing, MI 48824-1027\\

\noindent
E-mail:} ivanov@math.msu.edu

\end{document}